\documentclass[12pt]{article}
\usepackage{setspace}
\usepackage{amssymb,amsmath,amsthm}
\usepackage{graphicx}
\usepackage[english]{babel}
\newtheorem{Def}{Definition}
\newtheorem{Lem}{Lemma}
\newtheorem{Thm}[Lem]{Theorem}

\newtheorem{Ex}[Def]{Example}

\begin{document}

\begin{center}
\textbf{$(4,-(2n+5))$-TORUS KNOT WITH ONLY 1 NORMAL RULING}
\end{center}

\vspace{.01cm}
\begin{center}
\small{WATCHAREEPAN ATIPONRAT}
\end{center}

\vspace{.5cm}
\noindent
    A\scriptsize{BSTRACT.}\normalsize~ The main purpose of this paper is to provide an infinite family of counter examples of the open problem mentioned in \cite{Mep}. In particular, we present an infinite family of a particular Legendrian  $(4,-(2n+5))$-torus knot, for each $n \geq 0$, which has only 1 normal ruling, but do not satisfy the even number of clasps condition of Theorem 3 of $\cite{Mep}$. Thus, these normal rulings cannot imply the existence of a decomposable exact Lagrandian filling.
    
\vspace{1.3cm}
\begin{center}
\normalsize{1. I}\footnotesize{NTRODUCTION}
\end{center}

The \textbf{standard contact structure on $\mathbb{R}^3$} is a smooth 2-dimensional subbundle of the tangent bundle of $\mathbb{R}^3$ which is corresponding to ker $(dz-ydx)$. In this work, we consider smooth links in $\mathbb{R}^3$ which are everywhere tangent to the standard contact structure on $\mathbb{R}^3$, such links are called \textbf{Legendrian links}. If a Legendrian link has only 1 component, we may call it a Legendrian knot.

\begin{figure}
	\begin{center}
		\includegraphics[height=2in]{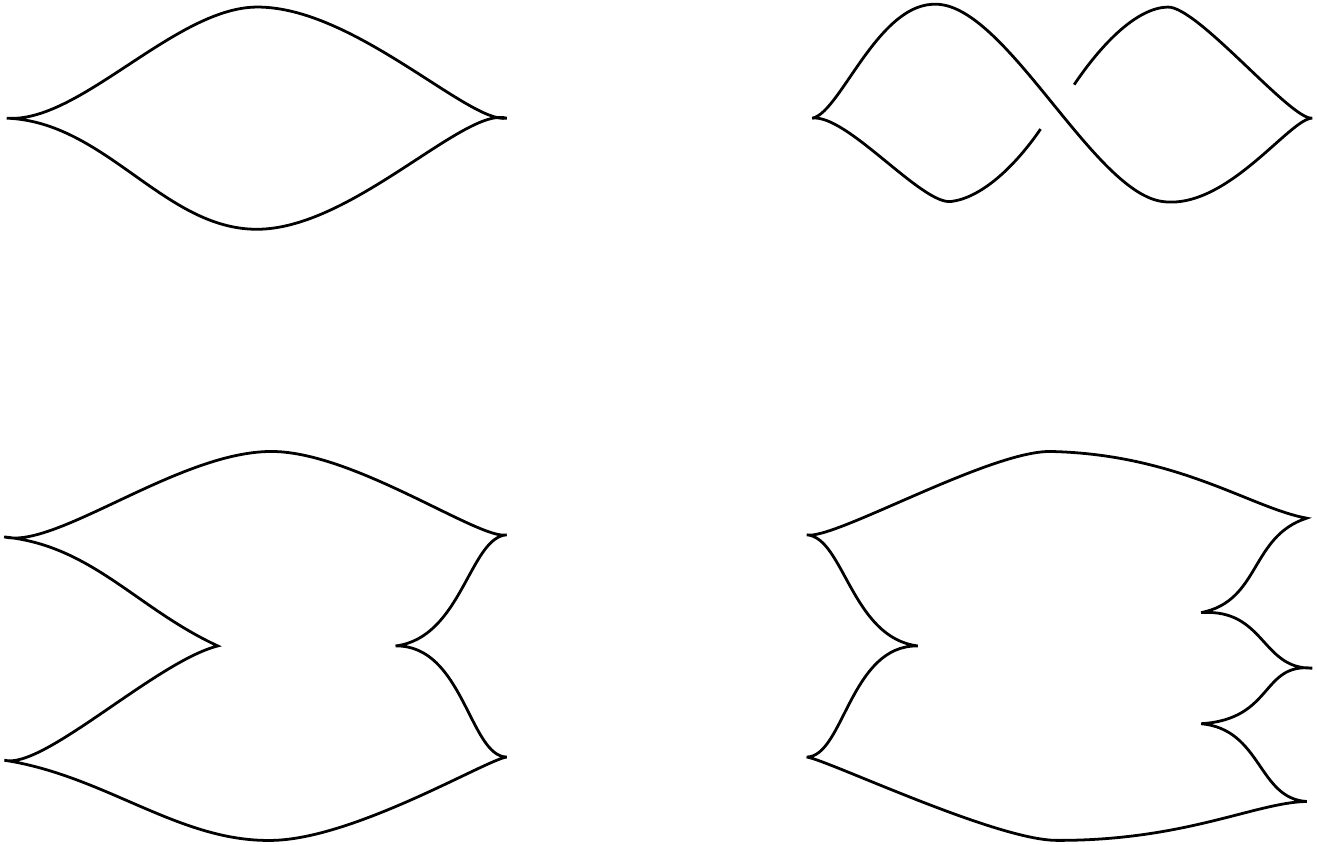}
	\end{center}
	\caption{Front projection of some Legendrian unknots.}
	\label{unknot}
	\vspace{.6cm}
\end{figure}

The \textbf{front projection} of a Legendrian link is the projection of the link into the $xz$-plane (here, we consider $\mathbb{R}^3$ with coordinate $(x,y,z)$). We will assume that the positive $y$-axis is pointing into the page so that every crossing of the front projection of a Legendrian link has the overpassing with less slope. Also, by small perturbation we may assume that front projection has only double points at self-intersections. It is proved in \cite{Geiges} that any self-intersection is transverse so the front projection has only finite number of self-intersections. In addition, a front projection of a Legendrian link has cusps instead of vertical tangencies. The front projection of some Legendrian unknots are illustrated in Figure \ref{unknot}.

Next, we introduce normal rulings. These are objects related to front diagram of Legendrian links. They are interesting since it has been shown in \cite{Chekanov Pushkar'} that the number of normal rulings of a Legendrian link is invariant under Legendrian isotopy. We will revisit this topic in Section 2.

On the other hand, an exact Lagrangian fillings of Legendrian links are particular surfaces in $\mathbb{R}^4$ with Legendrian links as their boundaries. The papers \cite{EHK} and \cite{Hayden} have suggested that there is a connection between the existence of normal rulings of Legendrian links and the existence of decomposable exact Lagrangian fillings of the links. As mentioned in the paragraph before Lemma 2 in \cite{Mep}, the existence of decomposable exact Lagrangian fillings is able to imply the existence of normal rulings. However, it is an open question if the occurrence of a normal ruling implies the existence of an exact (possibly non-orientable) Lagrangian filling. In this work, we would like to supply an infinite family of Legendrian links providing a negative answer to the question. In fact, we prove that there is an infinite family of a particular type of Legendrian knots which has only 1 normal ruling. To be precise, we have the following result at the end of section 2.

\begin{figure}
\begin{center}
\includegraphics{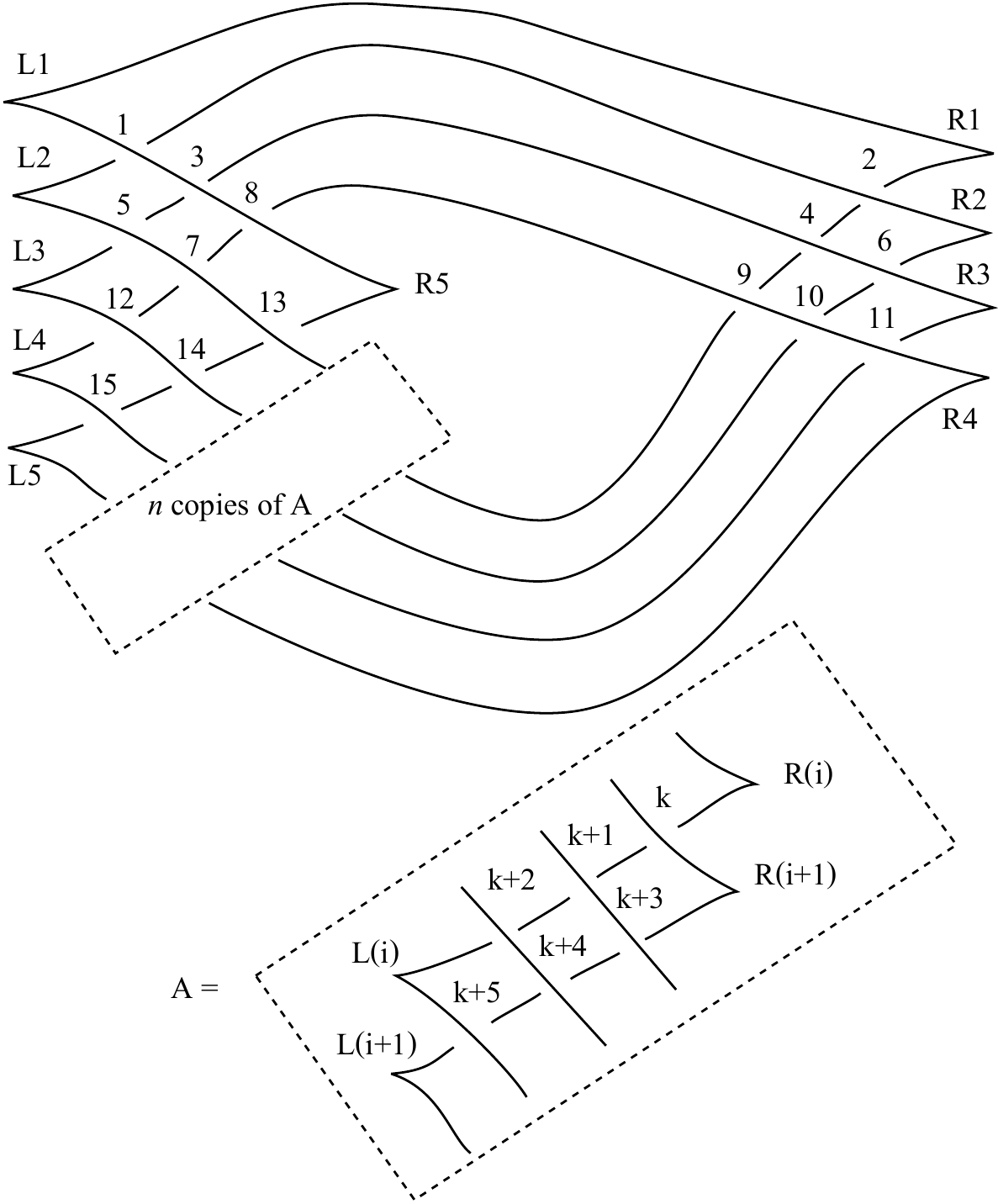}
\end{center}
\caption{The Legendrian $(4,-(2n+5))$-torus knot, $n \geq 0$.}
\label{tgen}
\end{figure}

\begin{Thm}
	\label{ngen}
	For any $n \geq 0$, Legendrian $(4,-(2n+5))$-torus knot, as in Figure \ref{tgen}, has only 1 normal ruling. Furthermore,  no member of this family of Legendrian knots has a decomposable exact Lagrandrian filling.
\end{Thm}

We note here that all results in this paper are coming from \cite{Met}. In addition, we will use this paper to provide an infinite family of counter examples in \cite{Mep}.

\vspace{0.3cm}
\noindent
\textbf{Acknowledgments.} The author would like to thank William Menasco for his support which makes this work possible. In addition, the author would like to thank Lenhard Ng for an introduction to this topic.

\begin{center}
\normalsize{2. N}\footnotesize{ORMAL RULINGS}
\end{center}

\begin{figure}[h]
\begin{center}
\vspace{1cm}
\includegraphics[width=2.5in]{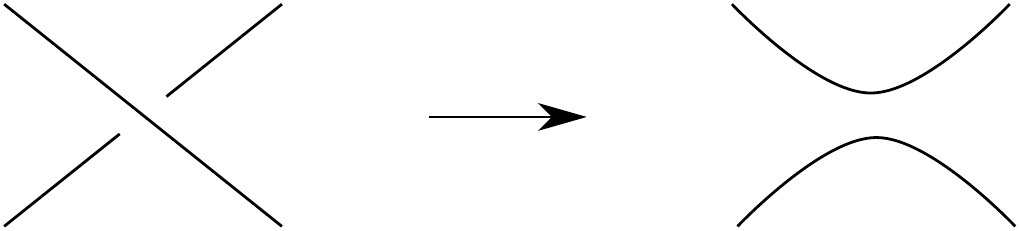}
\end{center}
\caption{Resolution at a crossing.}
\label{resolution}
\vspace{.3cm}
\end{figure}

Suppose we have a front diagram $K$ of a Legendrian link. By regular isotopy, we may assume from now on that its cusps and crossings have distinct $x$-coordinates. We consider a subset $\rho$ of the set of all crossings of $K$. Then we perform resolution, see Figure \ref{resolution}, at each crossing in $\rho$ so that we obtain a resulting front diagram $K'$. We call $\rho$ a \textbf{normal ruling} if the followings hold:

(1) each component of $K'$ has one left cusp, one right cusp and no self-intersections;

(2) horizontal strands at each resolution belong to different components in $K'$; and

\begin{figure}[h]
\begin{center}
\vspace{0.6cm}
\includegraphics[width=3.2in]{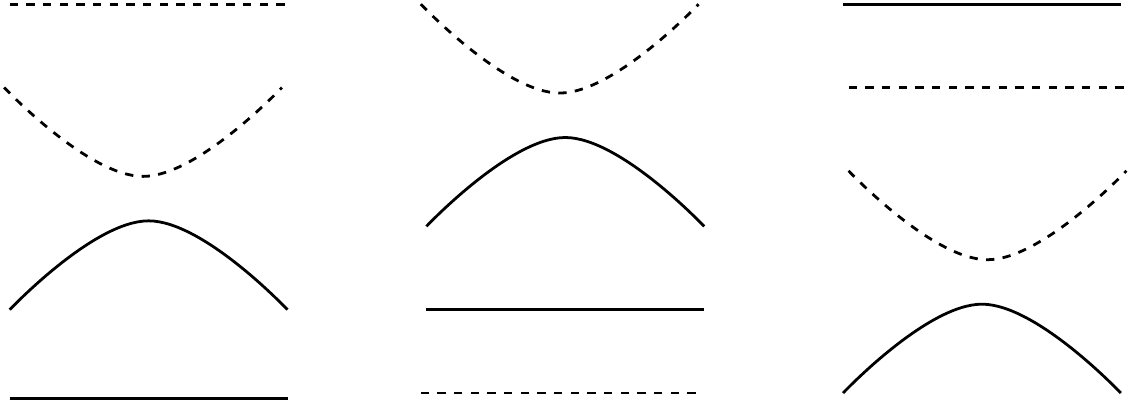}
\end{center}
\caption{Possible vertical slices at resolution when considering only two resulting components involved.}
\label{normality}
\vspace{.3cm}
\end{figure}

(3) in the vertical slice (constant $x$-coordinate) passing through each resolution, the two eyes meeting at the resolution must be one of the three cases in Figure \ref{normality}.

\begin{figure}
\begin{center}
\includegraphics{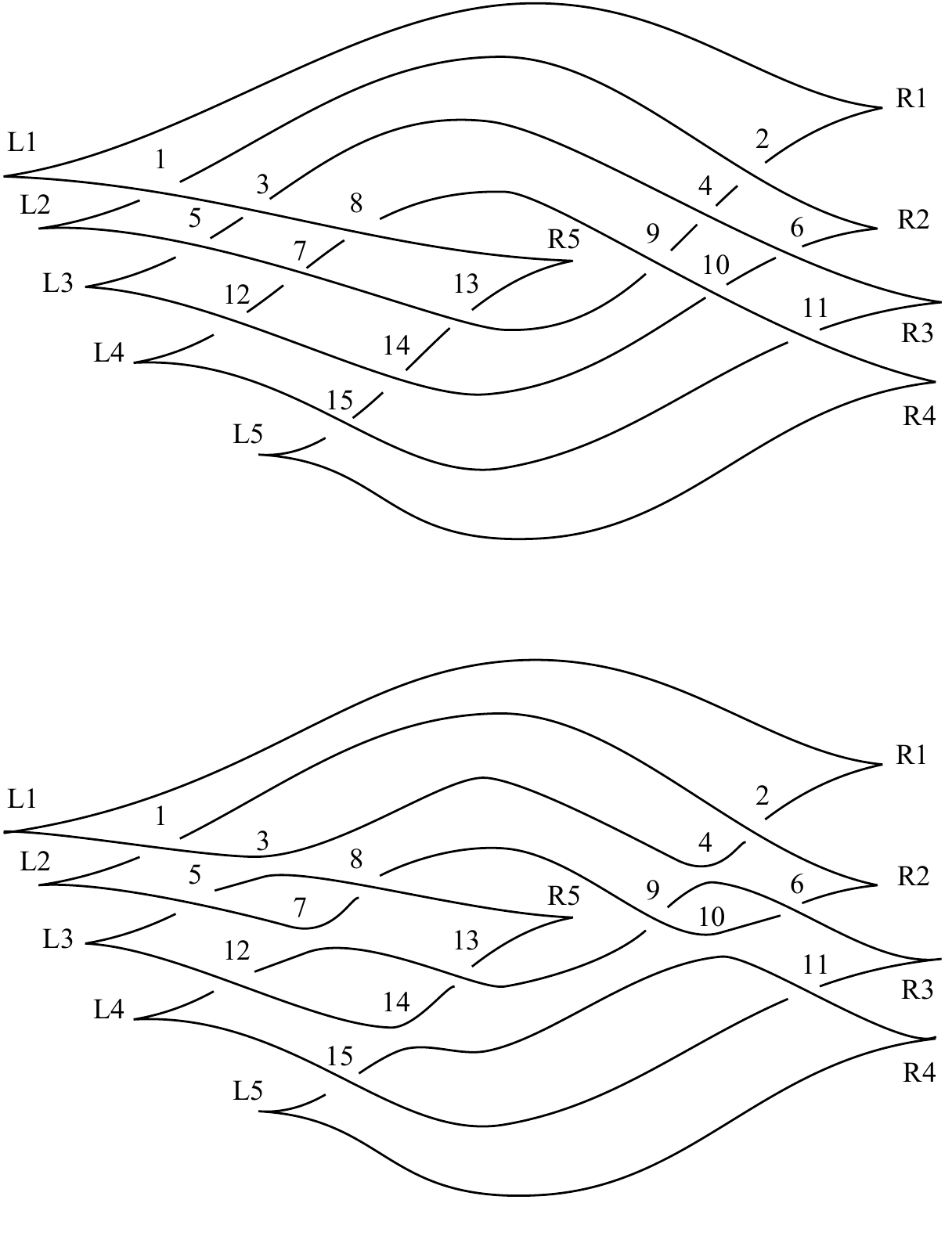}
\end{center}
\caption{The Legendrian $(4,-5)$-torus knot and the resolution of its only normal ruling.}
\label{m15}
\end{figure}

If $\rho$ is a normal ruling, then all crossings in $\rho$ are called \textbf{switches} and $K'$ is the \textbf{resolution of $\rho$} while each component of $K'$ is named an \textbf{eye}. Moreover, (3) is the \textbf{normality condition}, and we say a Legendrian link has a normal ruling if its front diagram admitting a normal ruling. Now, we give important examples of Legendrian knots with a normal ruling.
\begin{Ex}
\label{e1}
Legendrian $(4,-5)$-torus knot from Figure \ref{m15} has only 1 normal ruling.
\end{Ex}

\begin{figure}[h]
\begin{center}
\includegraphics[width=4.8in]{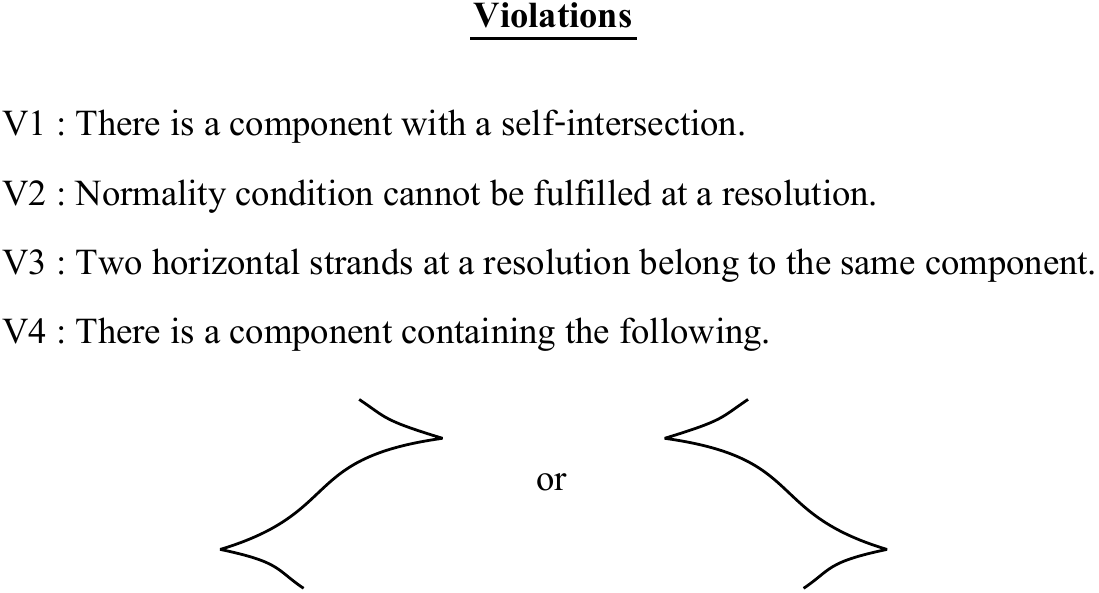}
\end{center}
\caption{Violations for being a normal ruling in Figure \ref{tab1} - \ref{tabc5}, \ref{tab4} - \ref{tabc7}, \ref{tabe1}, \ref{tabe2}, \ref{tabce}, \ref{tabo1}, \ref{tabo2}, \ref{tabco}.}
\label{vio}
\vspace{.3cm}
\end{figure}

\begin{figure}
\begin{center}
\includegraphics[width=3.5in]{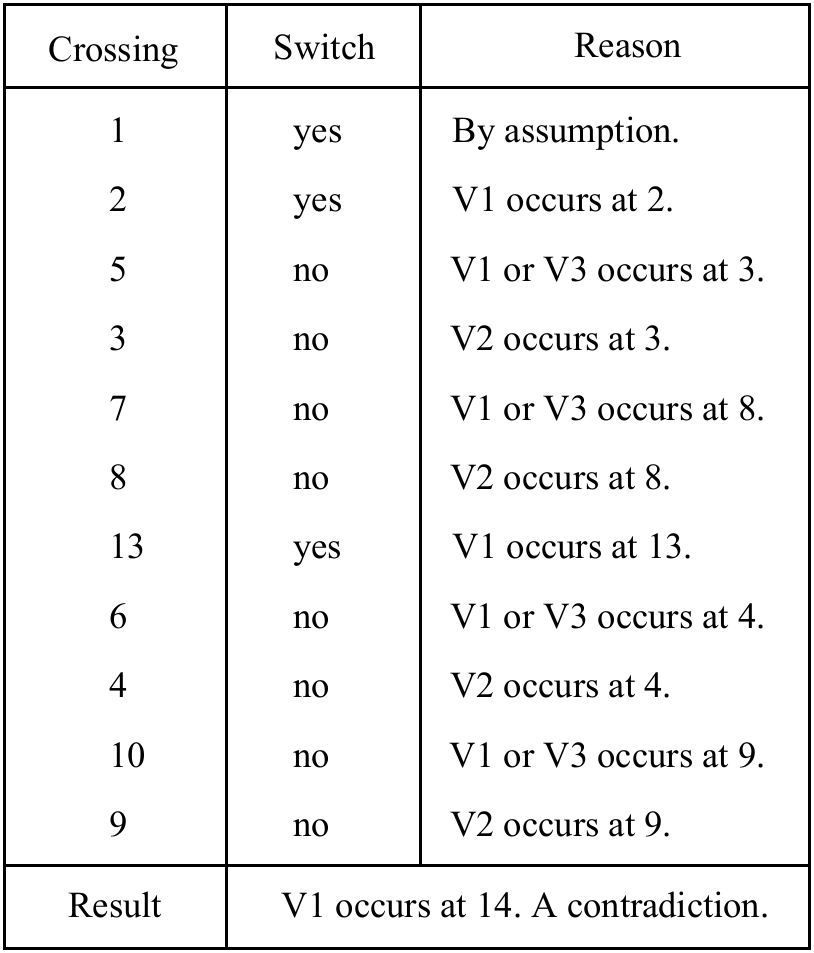}
\end{center}
\caption{Assuming 1 is a switch will give a contradiction.}
\label{tab1}
\end{figure}

\begin{figure}
\begin{center}
\includegraphics[width=3.5in]{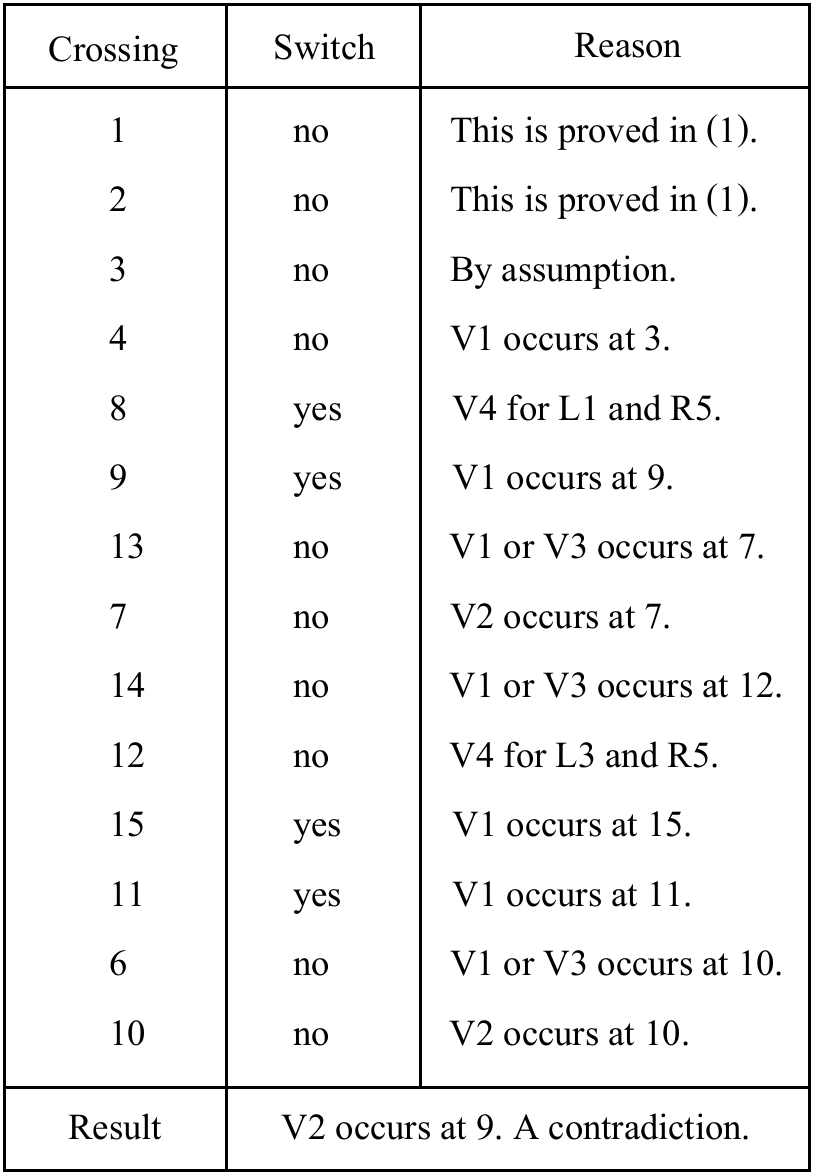}
\end{center}
\caption{Assuming 3 is not a switch will give a contradiction.}
\label{tab2}
\end{figure}

\begin{figure}
\begin{center}
\includegraphics[width=4.95in]{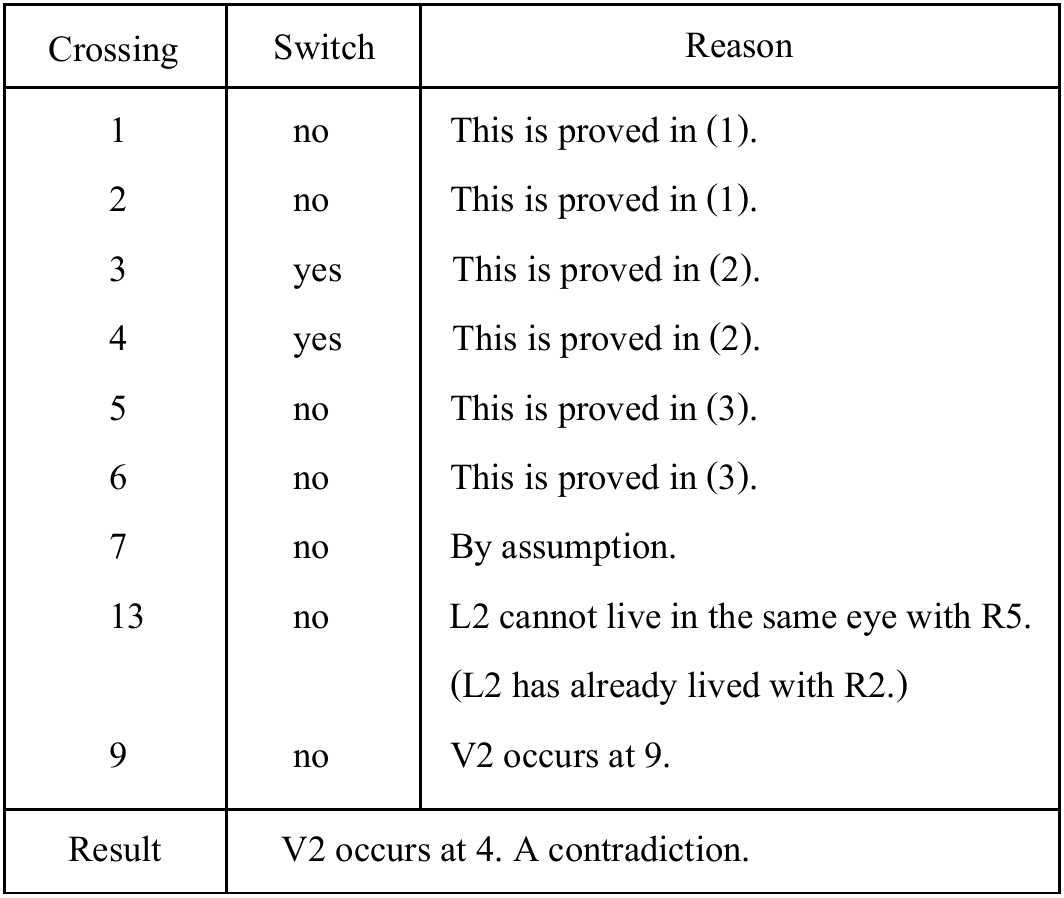}
\end{center}
\caption{Assuming 7 is not a switch will give a contradiction.}
\label{tab3}
\end{figure}

\begin{figure}
\begin{center}
\includegraphics[width=3.5in]{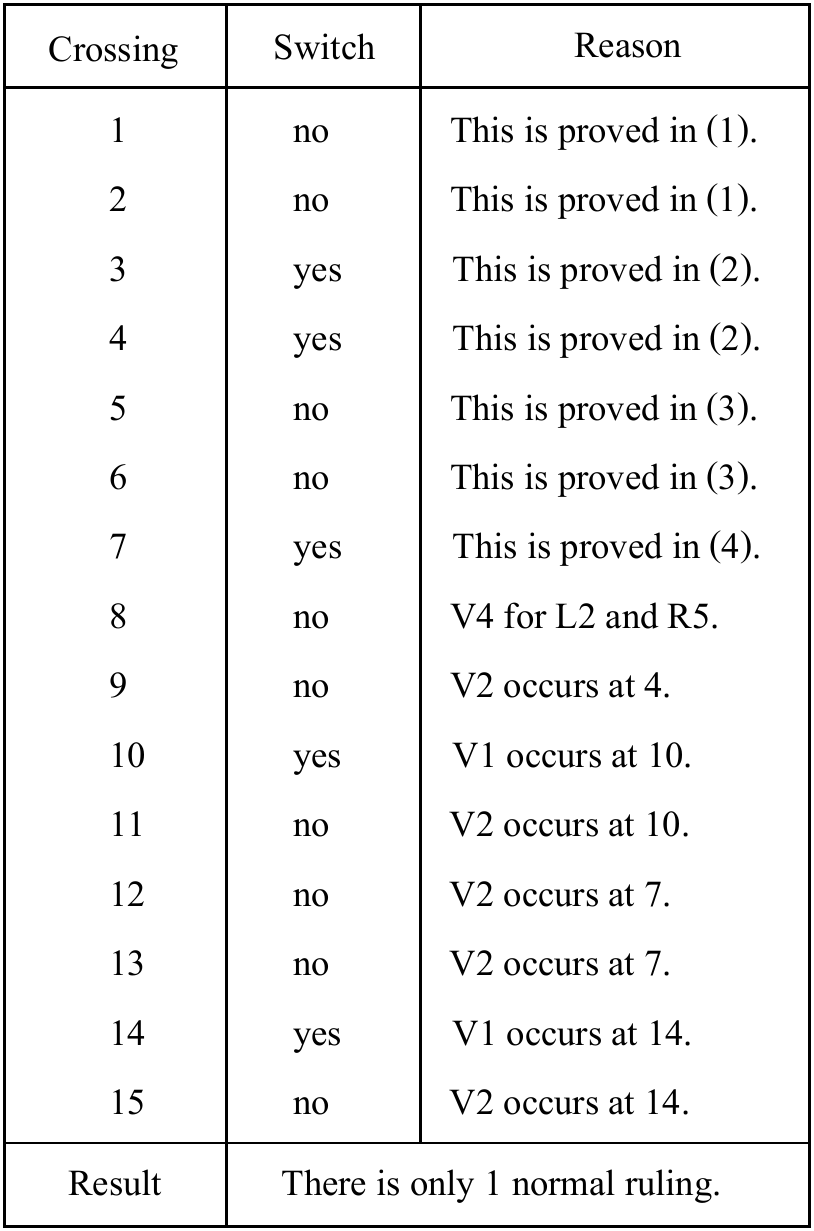}
\end{center}
\caption{There is only 1 normal ruling possible.}
\label{tabc5}
\end{figure}

\begin{proof}
First, it is not hard to check that \{3, 4, 7, 10, 14\} is a normal ruling as in the bottom of Figure \ref{m15}. Next, we show that it is the only possible normal ruling via some observations. Suppose we have a normal ruling. Then it cannot contain any of violations V1 - V4 for being a normal as shown in Figure \ref{vio}. It must satisfy the following.

(1) 1 and 2 are not switches.: Notice that we either have both 1 and 2 are switches or both are not switches. Suppose on the contrary that 1 is a switch. Then we will have its consequences and, at the end, a contradiction as illustrated in Figure \ref{tab1}. Also, because 1 and 2 are not switches, we have L2 and R2 live in the same eye.

(2) 3 and 4 are switches.: Since 1 and 2 are not switches, we either have both 3 and 4 are switches or both are not switches. Suppose on the contrary that 3 is not a switch. Then we will have its consequences and, at the end, a contradiction as illustrated in Figure \ref{tab2}.

(3) 5 and 6 are not switches.: If 5 or 6 is a switch, the normality condition fails at 3 or 4, which is impossible.

(4) 7 is a switch.: Suppose on the contrary that 7 is not a switch. Then we will have its consequences and, at the end, a contradiction as illustrated in Figure \ref{tab3}.

By (1) - (4), we only have 1 normal ruling possible as discussed in Figure \ref{tabc5}.

Thus there is exactly 1 normal ruling.
\end{proof}

\begin{figure}
\begin{center}
\includegraphics{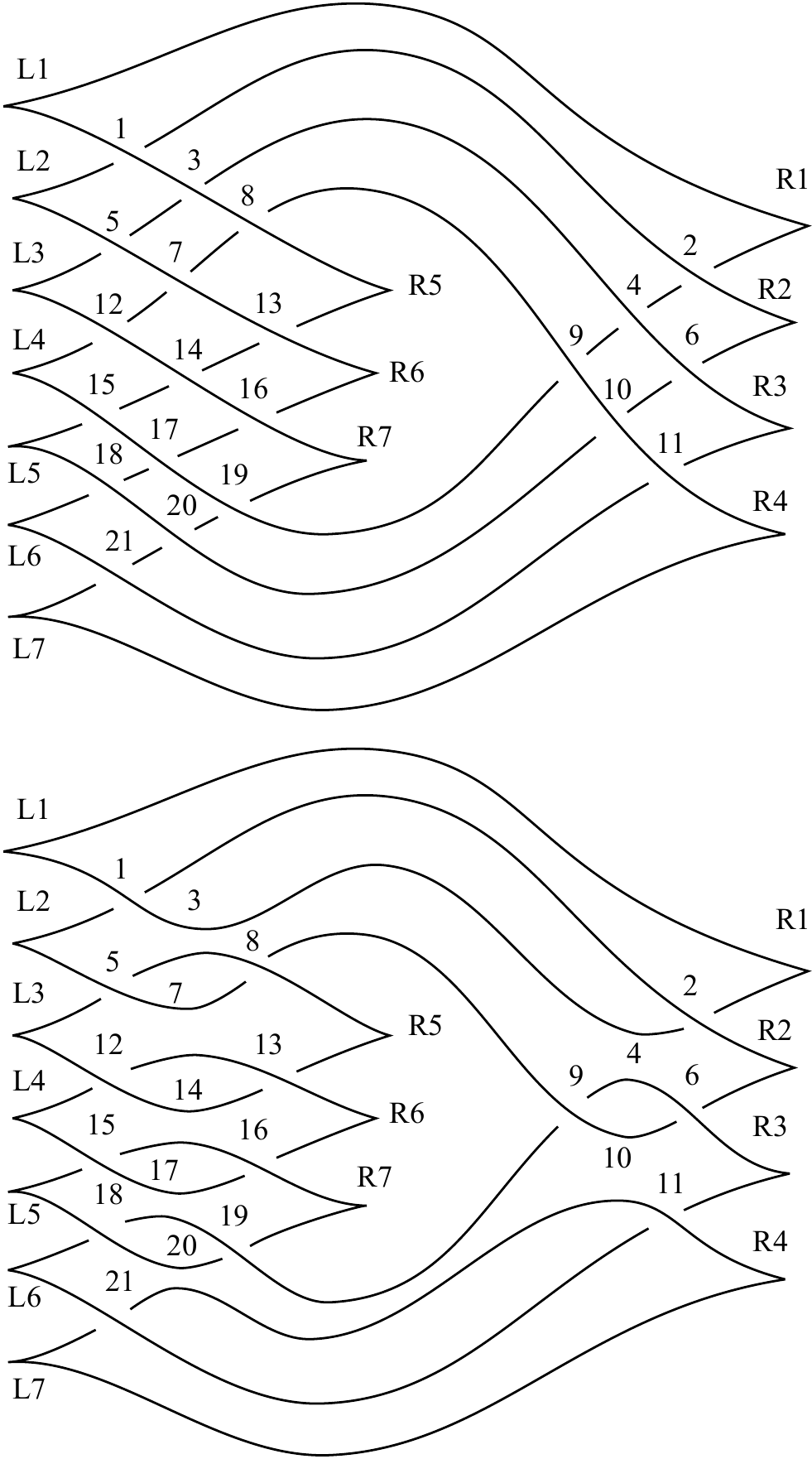}
\end{center}
\caption{The Legendrian $(4,-7)$-torus knot and the resolution of its only normal ruling.}
\label{m21}
\end{figure}

\newpage
\begin{Ex}
\label{e11}
Legendrian $(4,-7)$-torus knot from Figure \ref{m21} has only 1 normal ruling.
\end{Ex}

\begin{figure}
\begin{center}
\includegraphics[width=3.5in]{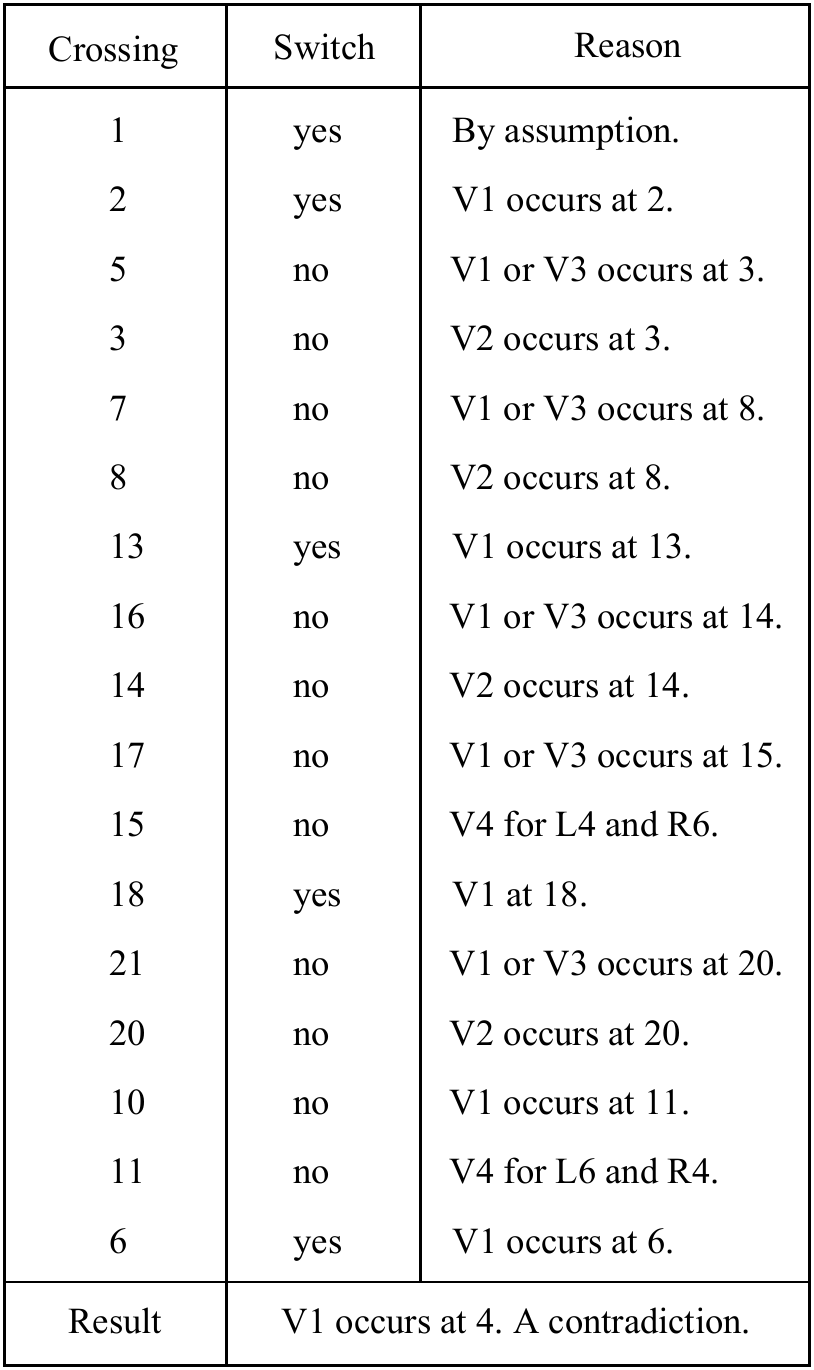}
\end{center}
\caption{Assuming 1 is a switch will give a contradiction.}
\label{tab4}
\end{figure}

\begin{figure}
\begin{center}
\includegraphics[width=3.5in]{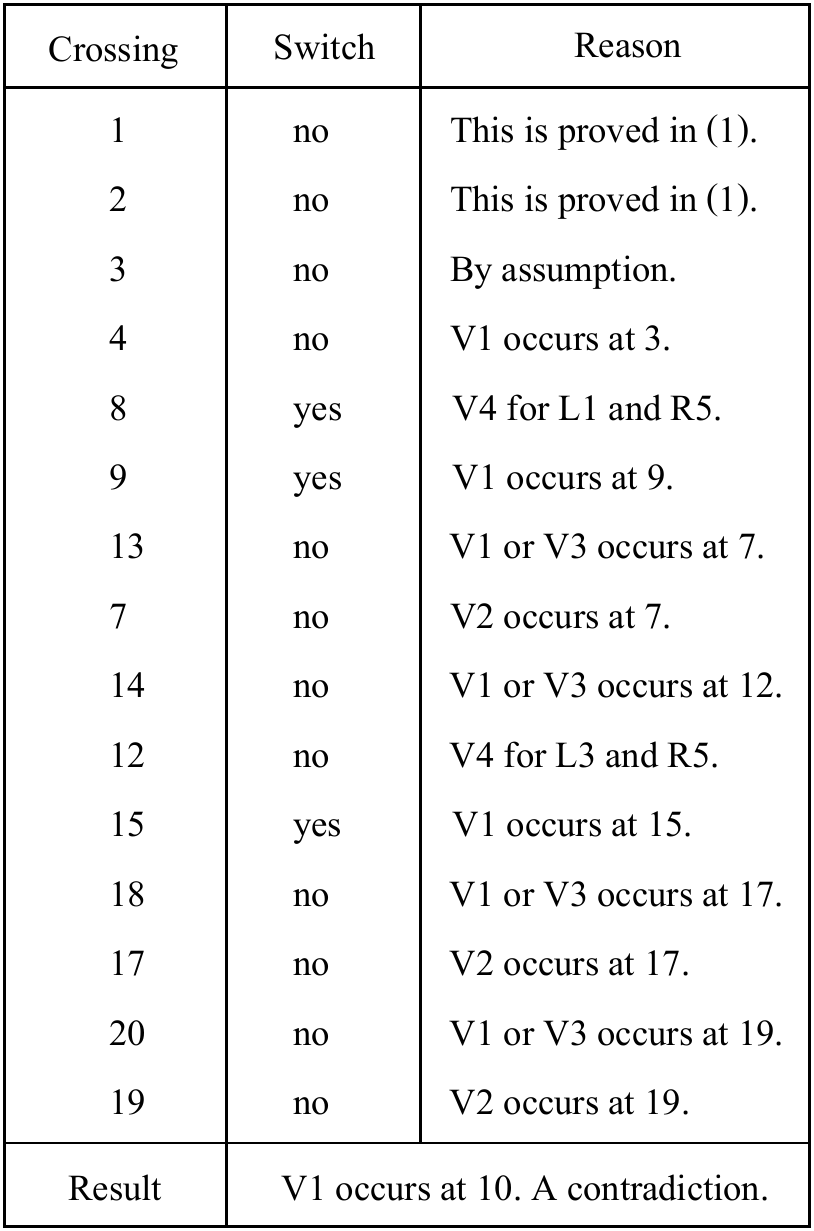}
\end{center}
\caption{Assuming 3 is not a switch will give a contradiction.}
\label{tab5}
\end{figure}

\begin{figure}
\begin{center}
\includegraphics[width=3.5in]{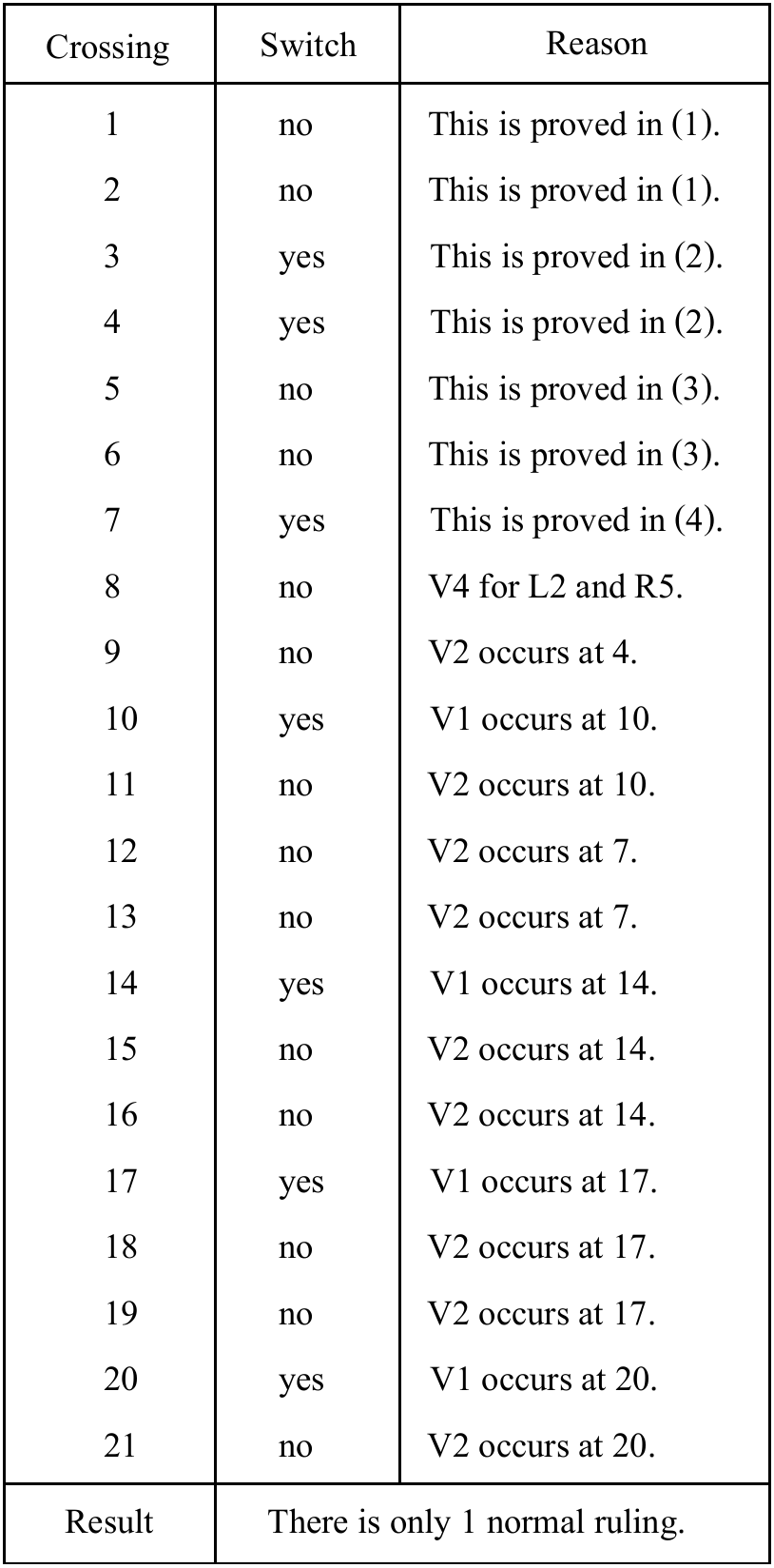}
\end{center}
\caption{There is only 1 normal ruling possible.}
\label{tabc7}
\end{figure}

\begin{proof}
First, it is not hard to check that \{3, 4, 7, 10, 14, 17, 20\} is a normal ruling as in the bottom of Figure \ref{m21}. Next, we show that it is the only possible normal ruling via some observations. Suppose we have a normal ruling. Then it must satisfy the following.

(1) 1 and 2 are not switches.: Notice that we either have both 1 and 2 are switches or both are not switches. Suppose on the contrary that 1 is a switch. Then we will have its consequences and, at the end, a contradiction as illustrated in Figure \ref{tab4}. Also, because 1 and 2 are not switches, we have L2 and R2 live in the same eye.

(2) 3 and 4 are switches.: Since 1 and 2 are not switches, we either have both 3 and 4 are switches or both are not switches. Suppose on the contrary that 3 is not a switch. Then we will have its consequences and, at the end, a contradiction as illustrated in Figure \ref{tab5}.

(3) 5 and 6 are not switches.: If 5 or 6 is a switch, the normality condition fails at 3 or 4, which is impossible.

(4) 7 is a switch.: Suppose on the contrary that 7 is not a switch. Then we will have its consequences and, at the end, a contradiction as illustrated in Figure \ref{tab3}.

By (1) - (4), we only have 1 normal ruling possible as discussed in Figure \ref{tabc7}.

Thus there is exactly 1 normal ruling.
\end{proof}

\begin{figure}
\begin{center}
\includegraphics[width=4.9in]{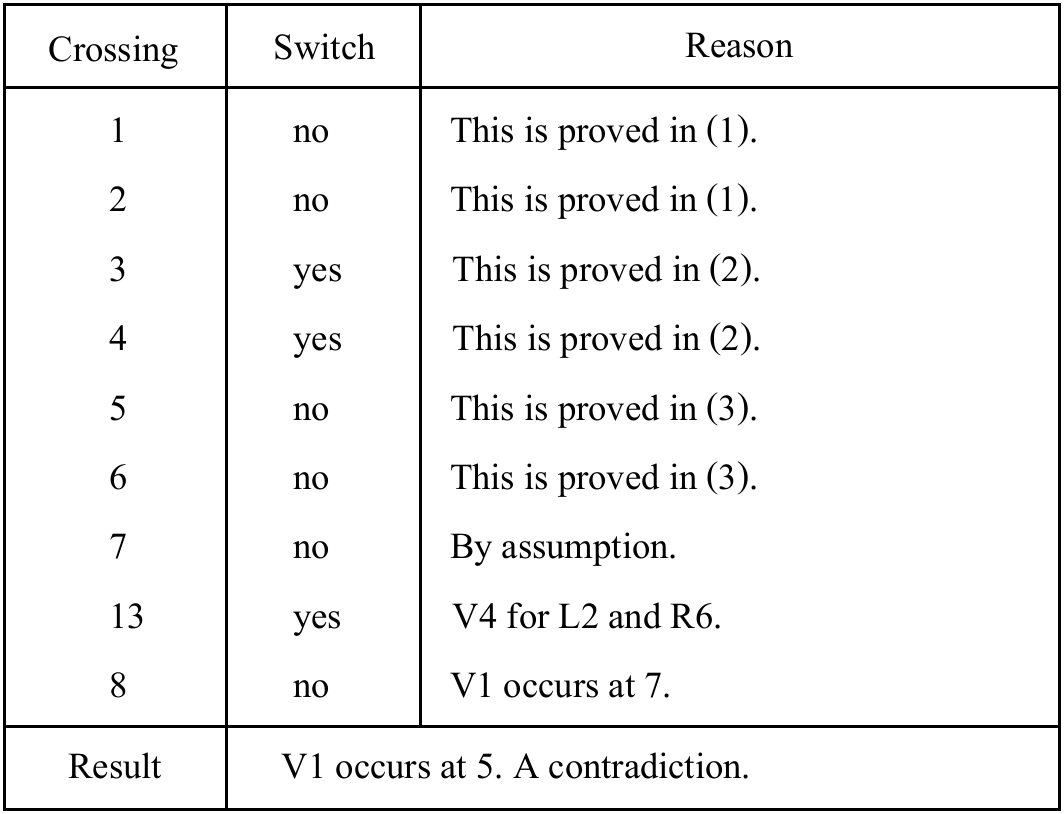}
\end{center}
\caption{Assuming 7 is not a switch will give a contradiction.}
\label{tab6}
\end{figure}

\begin{figure}
\begin{center}
\includegraphics{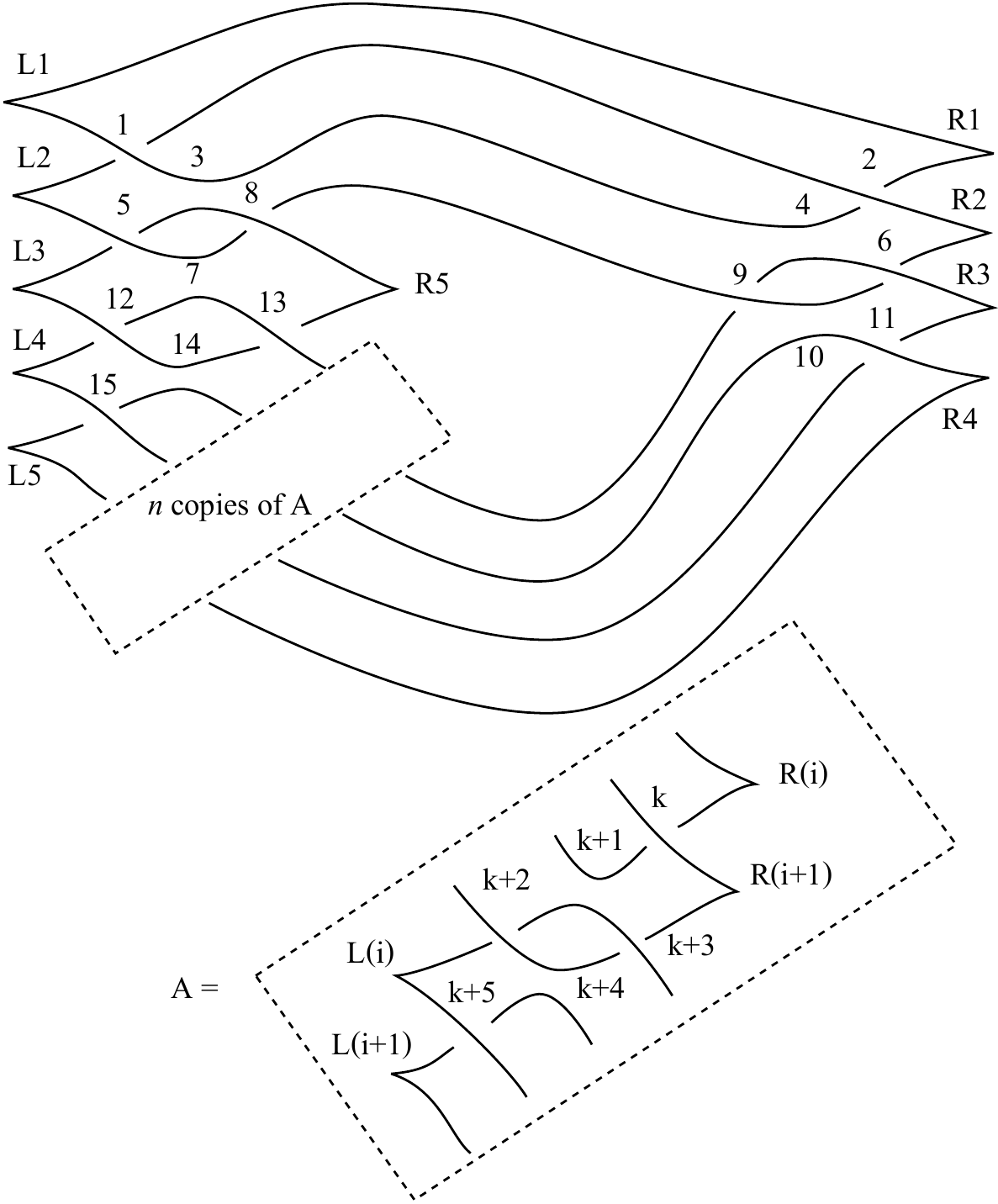}
\end{center}
\caption{The resolution of the only normal ruling of the Legendrian $(4,-(2n+5))$-torus knot, $n \geq 0$.}
\label{tgenr}
\end{figure}

\begin{figure}
\begin{center}
\includegraphics[width=4.6in]{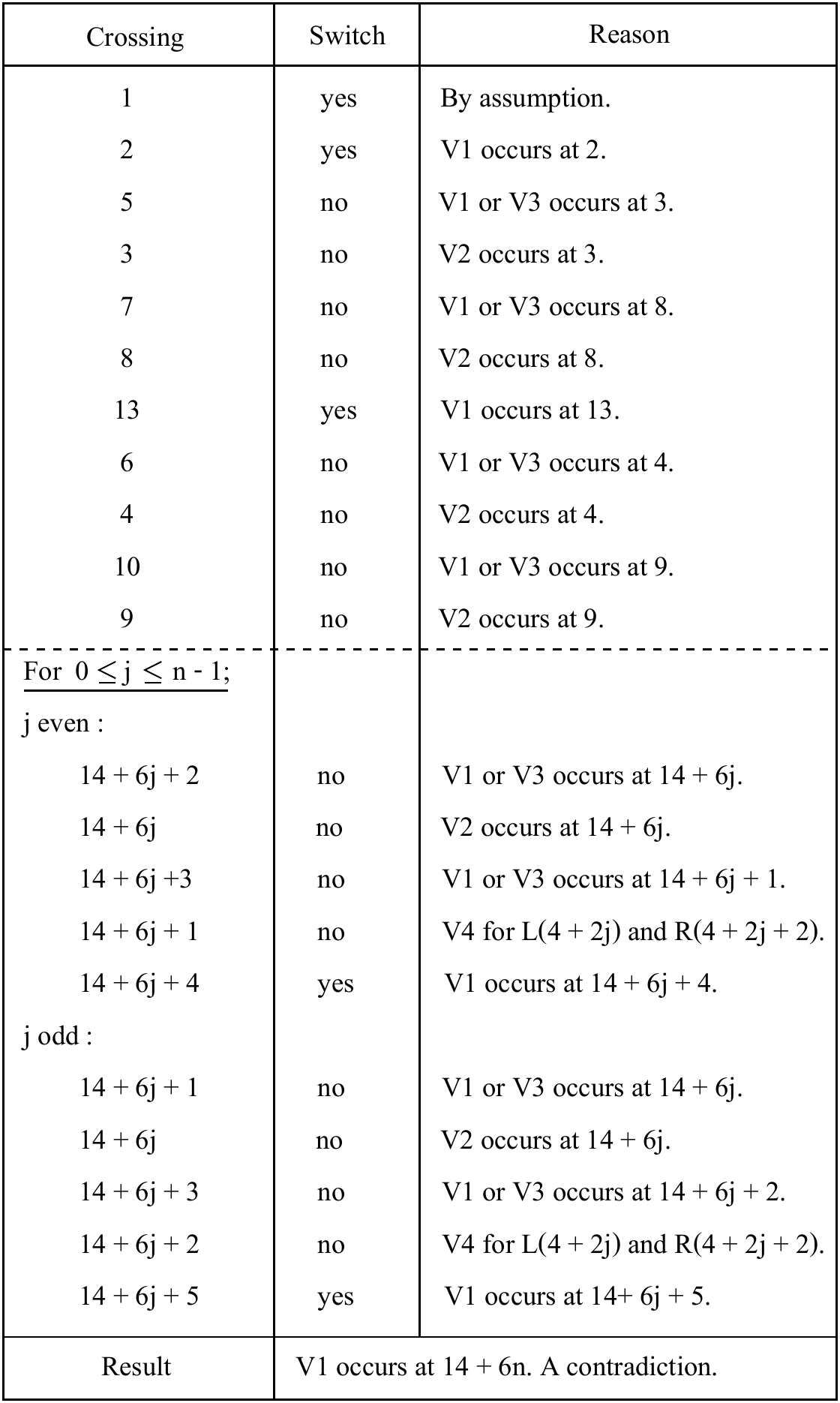}
\end{center}
\caption{Assuming 1 is a switch when $n$ is even will give a contradiction.}
\label{tabe1}
\end{figure}

\begin{figure}
\begin{center}
\includegraphics[width=4.4in]{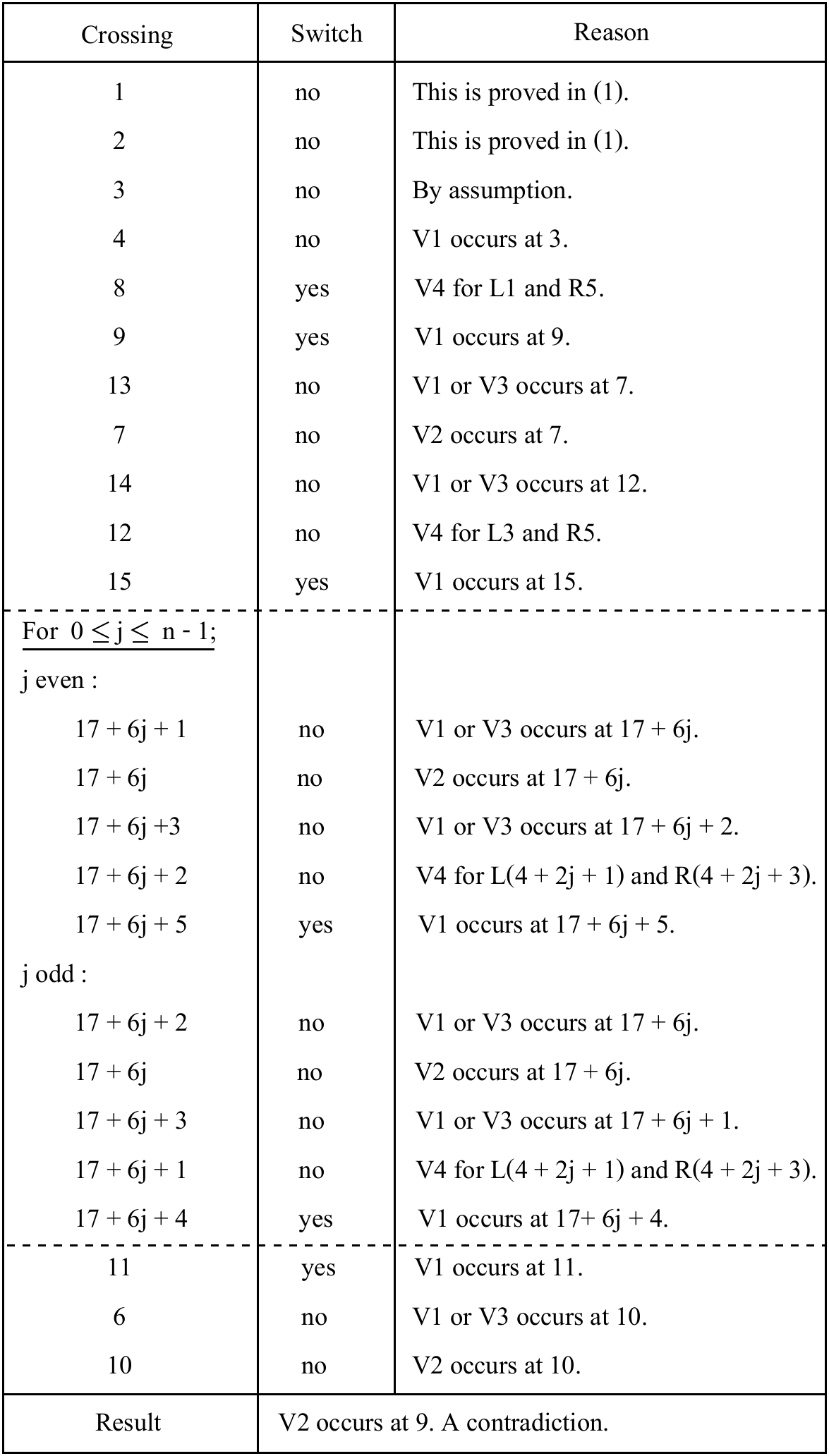}
\end{center}
\caption{Assuming 3 is not a switch when $n$ is even will give a contradiction.}
\label{tabe2}
\end{figure}

\begin{figure}
\begin{center}
\includegraphics[width=4.5in]{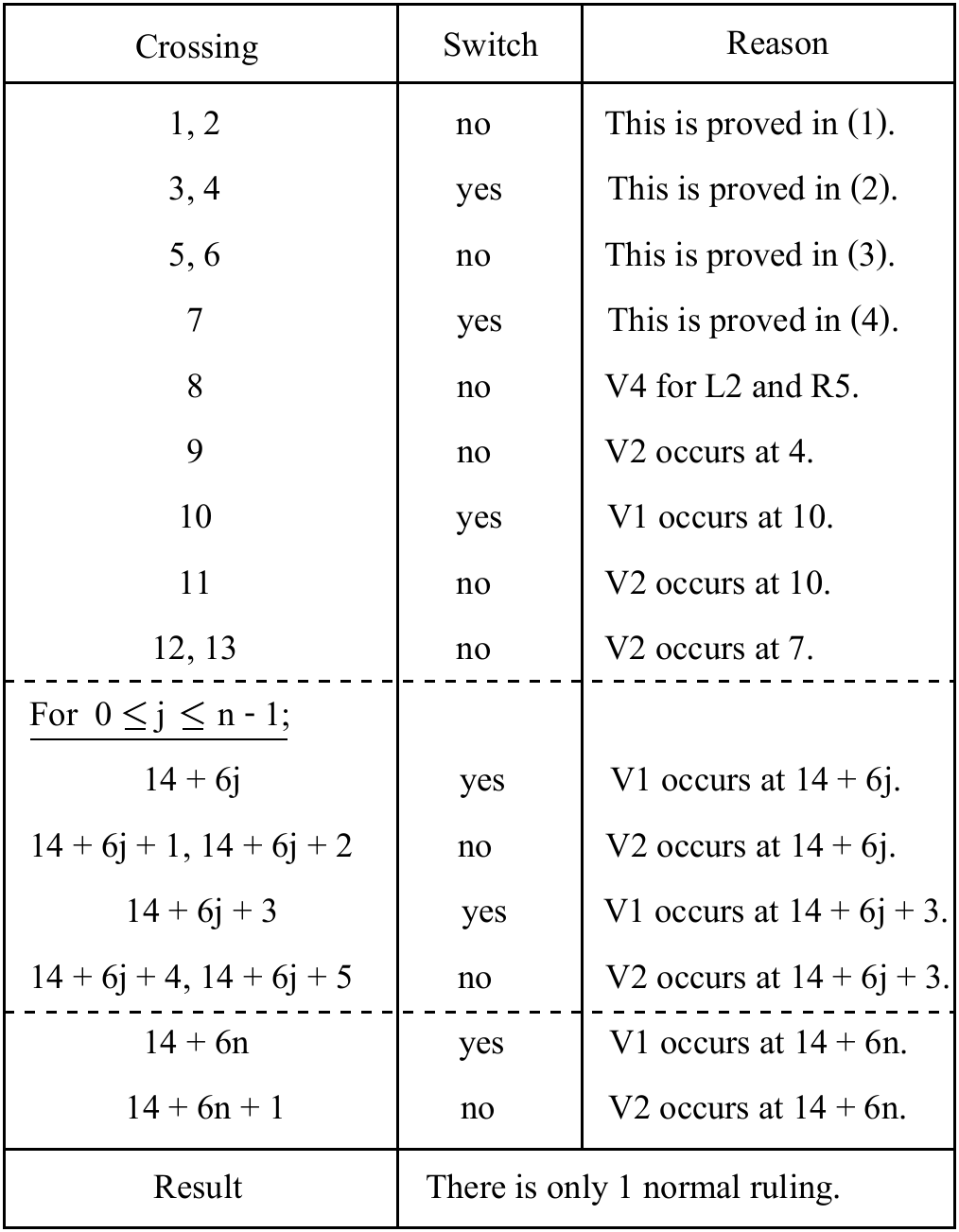}
\end{center}
\caption{There is only 1 normal ruling possible when $n$ is even.}
\label{tabce}
\end{figure}

\begin{figure}
\begin{center}
\includegraphics[width=4.5in]{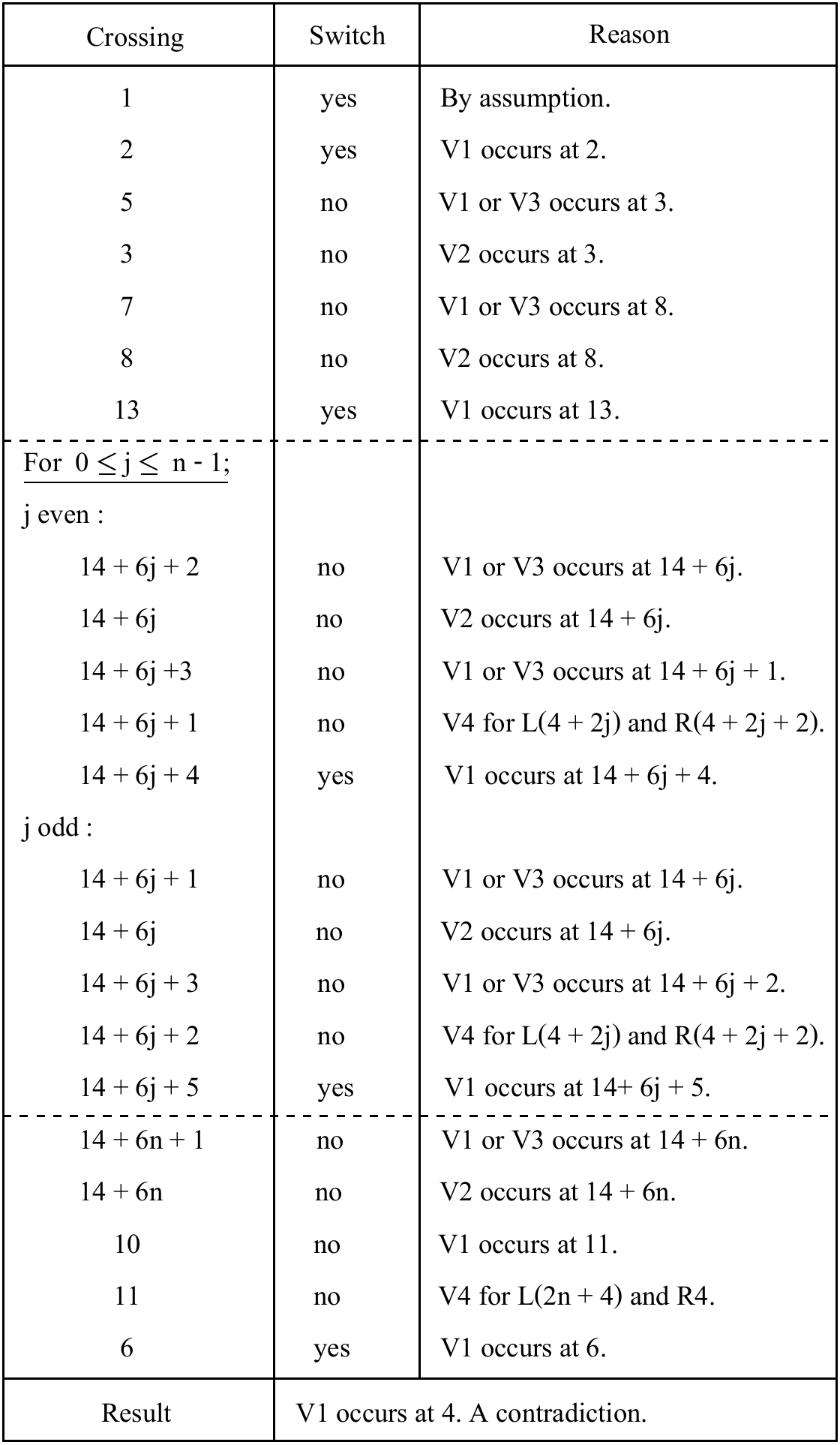}
\end{center}
\caption{Assuming 1 is a switch when $n$ is odd will give a contradiction.}
\label{tabo1}
\end{figure}

\begin{figure}
\begin{center}
\includegraphics[width=4.5in]{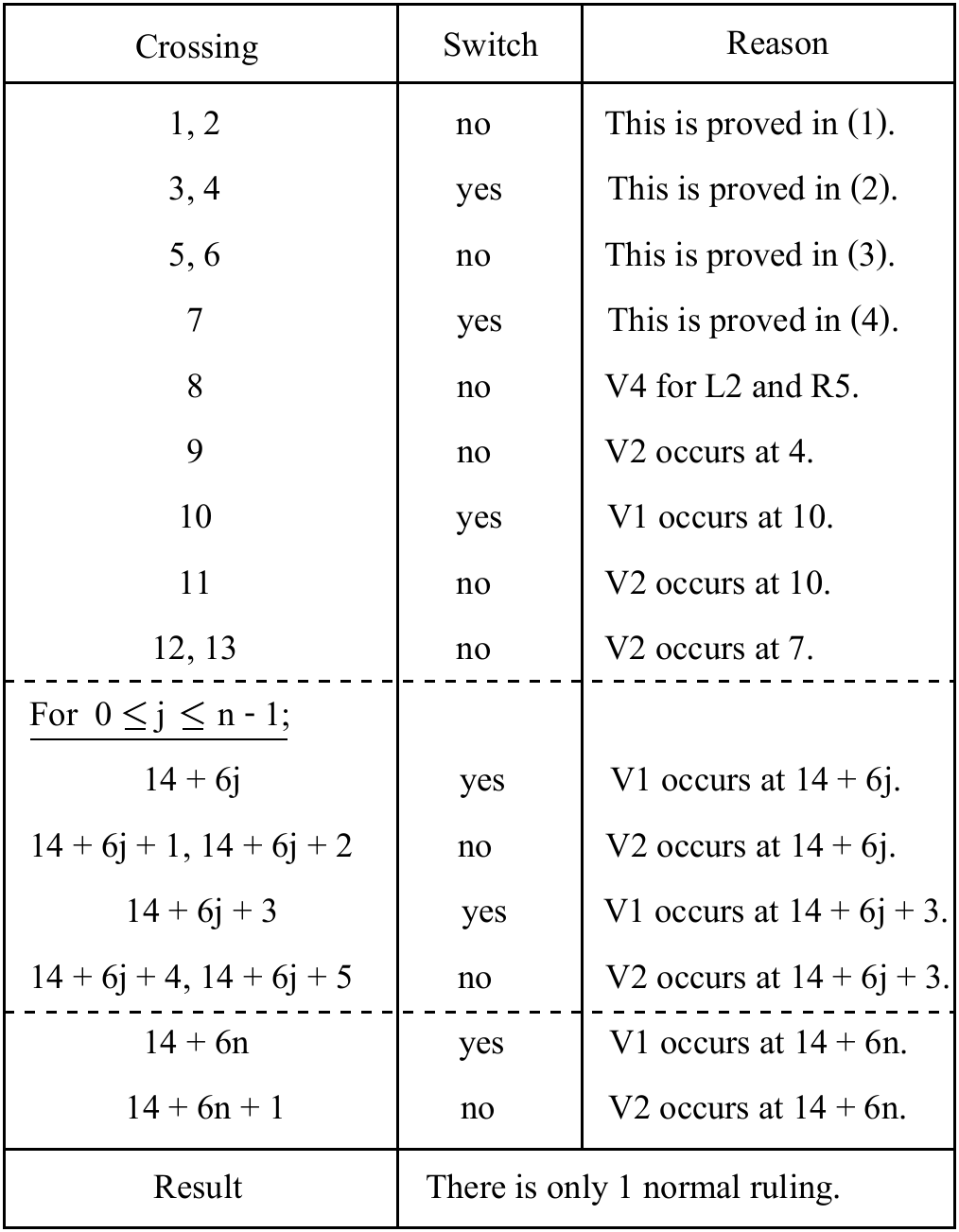}
\end{center}
\caption{There is only 1 normal ruling possible when $n$ is odd.}
\label{tabco}
\end{figure}

\newpage
We may generalize Example \ref{e1} - \ref{e11} to obtain Theorem 1 as follows.
\vspace{0.2cm}
\begin{proof}[Proof of Theorem 1.]
We prove by induction on $n$. The case $n =$ 0, 1 are shown in previous examples. Now, suppose the statement is true when $0 \leq n \leq N$ for some $N \geq 1$. We want to prove the statement is true when $n = N+1$. First, \{3, 4, 7, 10, 14, 14+3, 14+6, 14+9, ..., 14+$6n$-3, 14+$6n$\} is a normal ruling as in Figure \ref{tgenr}. Next, we show that it is the only possible normal ruling. We have 2 cases to consider.
~~
\\
\underline{$n = N+1$ is even.}: Any normal ruling must satisfy the following.

(1) 1 and 2 are not switches.: Notice that we either have both 1 and 2 are switches or both are not switches. Suppose on the contrary that 1 is a switch. Then we will have its consequences and, at the end, a contradiction as illustrated in Figure \ref{tabe1}. Also, because 1 and 2 are not switches, we have L2 and R2 live in the same eye.

(2) 3 and 4 are switches.: Since 1 and 2 are not switches, we either have both 3 and 4 are switches or both are not switches. Suppose on the contrary that 3 is not a switch. Then we will have its consequences and, at the end, a contradiction as illustrated in Figure \ref{tabe2}.

(3) 5 and 6 are not switches.: If 5 or 6 is a switch, the normality condition fails at 3 or 4, which is impossible.

(4) 7 is a switch.: Since $n \geq 2$, we may use the prove of (4) from Example \ref{e11}. Suppose on the contrary that 7 is not a switch. Then we will have its consequences and, at the end, a contradiction as illustrated in Figure \ref{tab3}.

By (1) - (4), we only have 1 normal ruling possible as discussed in Figure \ref{tabce}.
~~
\\
\underline{$n = N+1$ is odd.}: Any normal ruling must satisfy the following.

(1) 1 and 2 are not switches.: Notice that we either have both 1 and 2 are switches or both are not switches. Suppose on the contrary that 1 is a switch. Then we will have its consequences and, at the end, a contradiction as illustrated in Figure \ref{tabo1}. Also, because 1 and 2 are not switches, we have L2 and R2 live in the same eye.

(2) 3 and 4 are switches.: Since 1 and 2 are not switches, we either have both 3 and 4 are switches or both are not switches. Suppose on the contrary that 3 is not a switch. Then we will have its consequences and, at the end, a contradiction as illustrated in Figure \ref{tabe2}.

(3) 5 and 6 are not switches.: If 5 or 6 is a switch, the normality condition fails at 3 or 4, which is impossible.

(4) 7 is a switch.: Since $n \geq 2$, we may use the prove of (4) from Example \ref{e11}. Suppose on the contrary that 7 is not a switch. Then we will have its consequences and, at the end, a contradiction as illustrated in Figure \ref{tab3}.

By (1) - (4), we only have 1 normal ruling possible as discussed in Figure \ref{tabco}.

Thus there is exactly 1 normal ruling.

Finally, the second part of this theorem is proved in Theorem 4 of $\cite{Mep}$.
\end{proof}

\begin{figure}
\begin{center}
\includegraphics[width=4.6in]{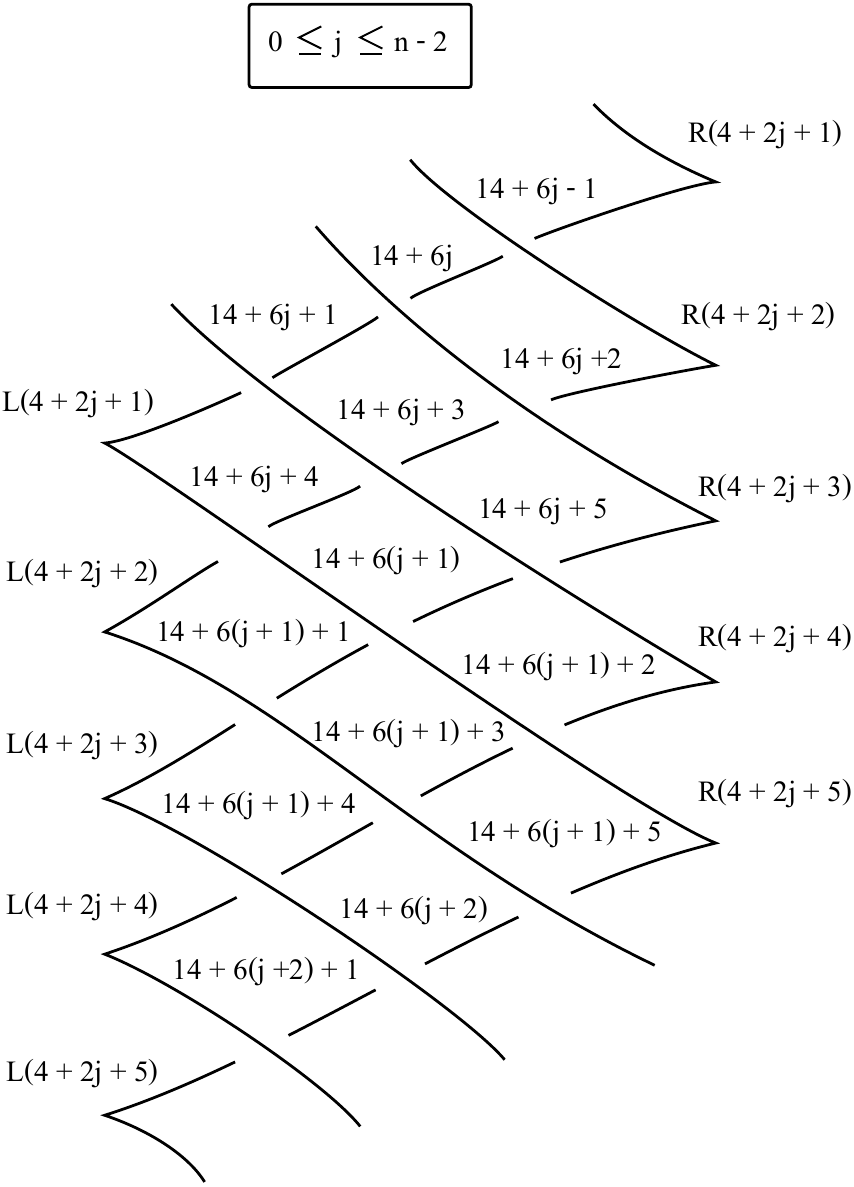}
\end{center}
\caption{Labeling crossings and cusps for tables in Figure \ref{tabe1} and \ref{tabce}.}
\label{tabe1c}
\end{figure}

\begin{figure}
\begin{center}
\includegraphics[width=4.6in]{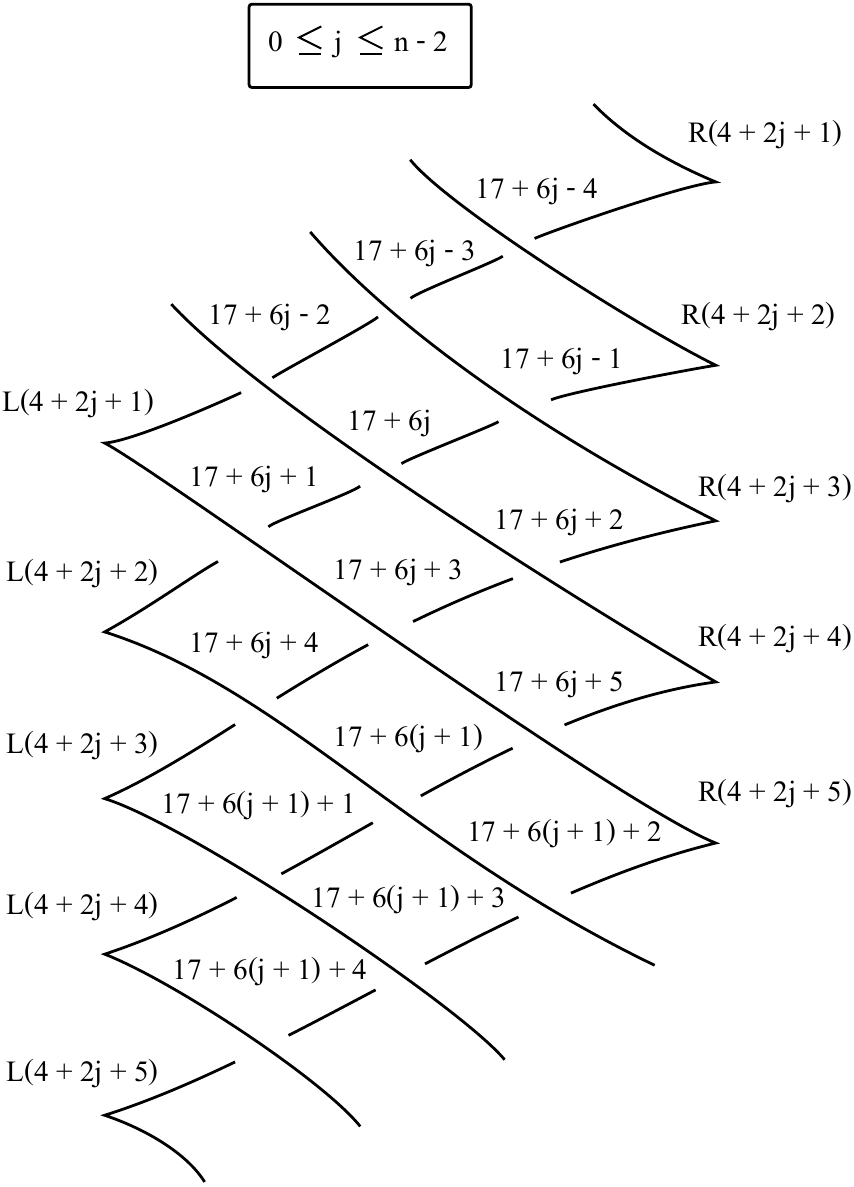}
\end{center}
\caption{Labeling crossings and cusps for table in Figure \ref{tabe2}.}
\label{tabe2c}
\end{figure}

\begin{figure}
\begin{center}
\includegraphics[width=4.6in]{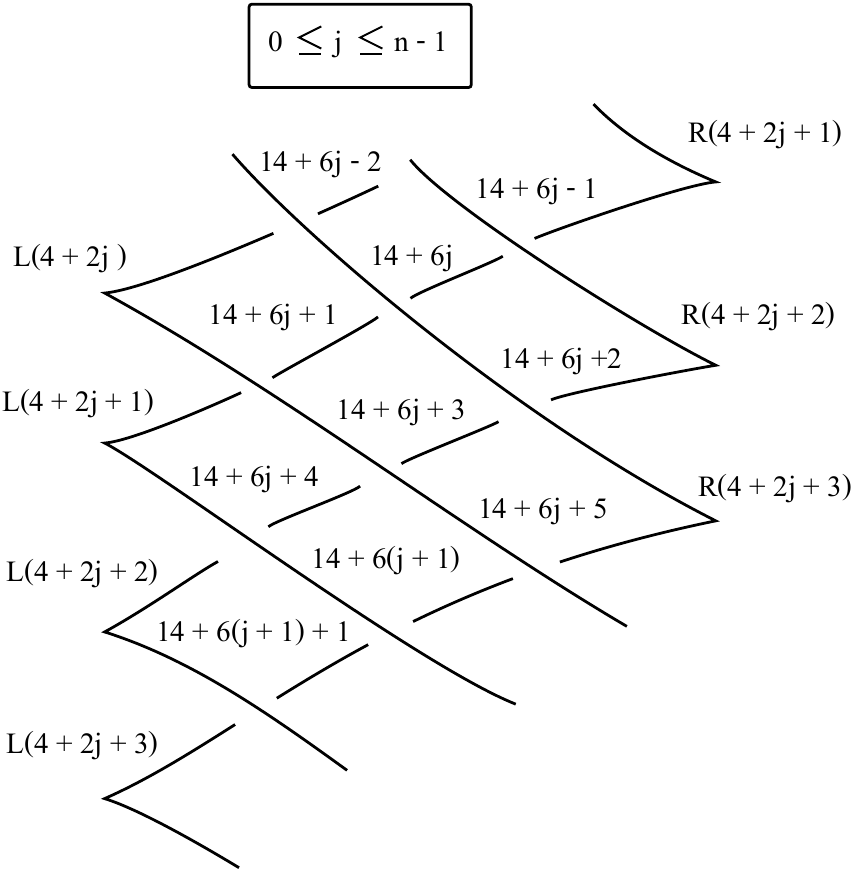}
\end{center}
\caption{Labeling crossings and cusps for table in Figure \ref{tabo1} and \ref{tabco}.}
\label{tabo1c}
\end{figure}

\begin{figure}
\begin{center}
\includegraphics[width=4.2in]{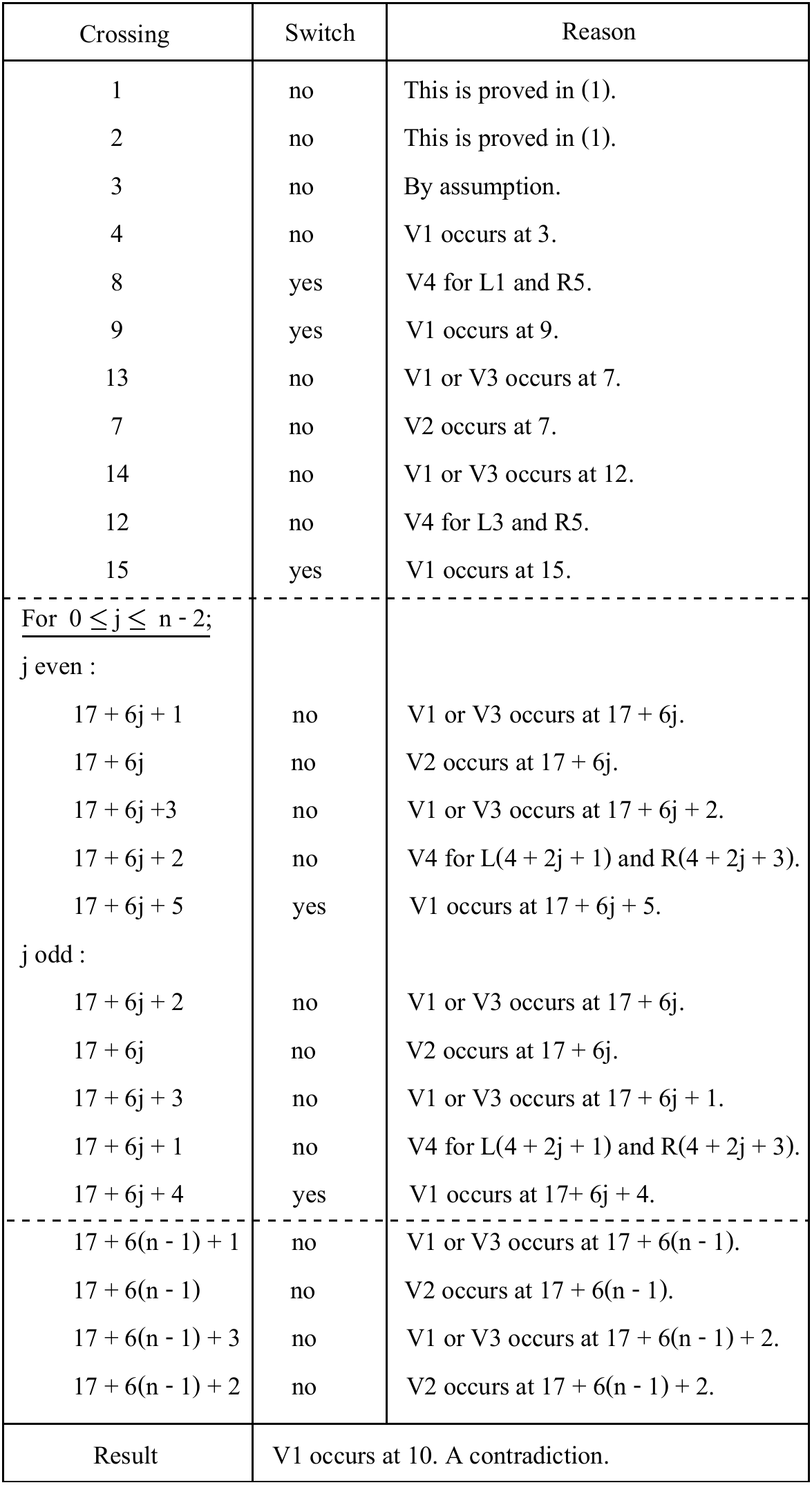}
\end{center}
\caption{Assuming 3 is not a switch when $n$ is odd will give a contradiction.}
\label{tabo2}
\end{figure}

\begin{figure}
\begin{center}
\includegraphics[width=4.6in]{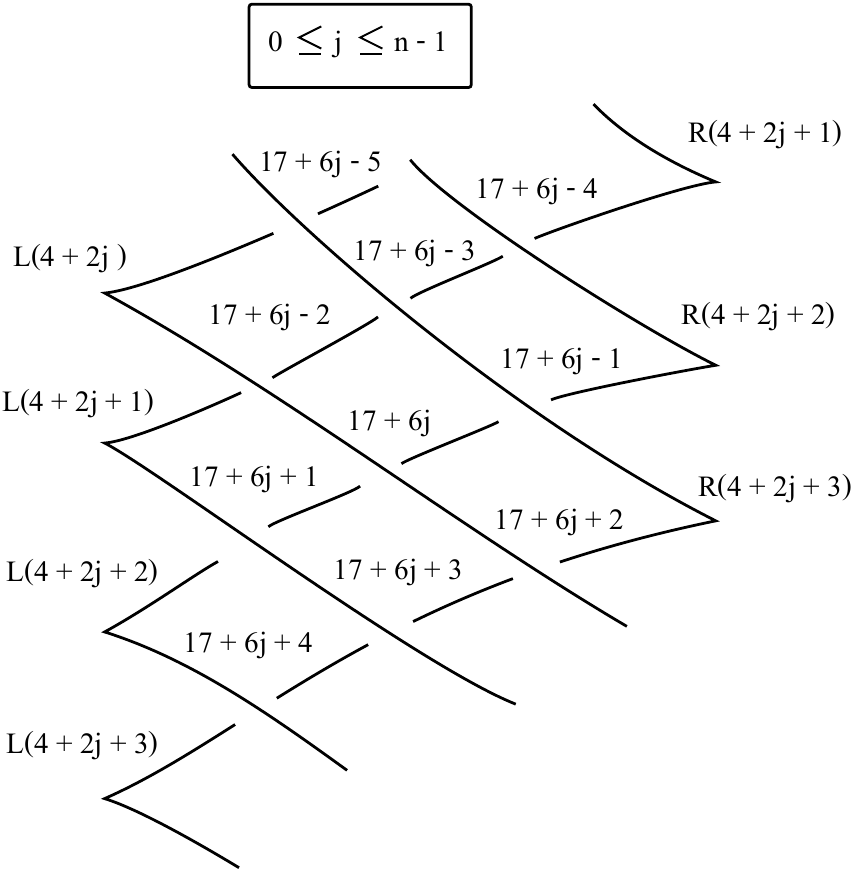}
\end{center}
\caption{Labeling crossings and cusps for table in Figure \ref{tabo2}.}
\label{tabo2c}
\end{figure}

\noindent 
Mathematics Department, Chiang Mai University\\
\textit{E-mail address}: atiponrat@hotmail.com

\end{document}